# Coset enumeration strategies

George Havas [*]
Key Centre for Software Technology
Department of Computer Science
University of Queensland
Queensland 4072
Australia

## Abstract

A primary reference on computer implementation of coset enumeration procedures is a 1973 paper of Cannon, Dimino, Havas and Watson. Programs and techniques described there are updated in this paper. Improved coset definition strategies, space saving techniques and advice for obtaining improved performance are included. New coset definition strategies for Felsch-type methods give substantial reductions in total cosets defined for some pathological enumerations. Significant time savings are achieved for coset enumeration procedures in general. Statistics on performance are presented, both in terms of time and in terms of maximum and total cosets defined for selected enumerations.

## 1 Introduction

Coset enumeration programs implement systematic procedures for enumerating the cosets of a subgroup $H$ of finite index in a group $G$, given a set of defining relations for $G$ and words generating $H$. Coset enumeration is the basis of important procedures for investigating finitely presented groups. Computer implementations are based on methods initially described by Todd and Coxeter [16]. A number of computer programs for coset enumeration have been described, including that by Cannon, Dimino, Havas and Watson [4].

Cannon *et al* include comprehensive references to earlier implementations, while later descriptions may be found in Neubüser [13], Leech [11] and Sims [15]. Neubüser and Leech both provide useful introductions to coset enumeration while Sims, in a substantial chapter of his forthcoming book, gives a formal account of coset enumeration in terms of automata and proves interesting new results on coset enumeration behaviour. This paper provides information on improved coset definition strategies and on space saving techniques and advice for obtaining enhanced performance from coset enumeration programs. More details on the experimental work which led to the selection of these methods is given in a technical report by Havas and Lian [7].

The procedures described in [4] remain a primary reference for coset enumeration program implementations including Felsch-type, Haselgrove-Leech-Trotter (HLT-type), and lookahead-type methods. Coset enumeration strategies are also discussed in [4]. Later work on lookahead methods is presented by Leon [12], while Sims [15] describes ten different strategies including both HLT-type and Felsch-type methods. This paper complements the valuable results in those with general advice on performance enhancement and specific strategies for Felsch-type methods.

In the group theory language Cayley (see Cannon [3]), the Todd-Coxeter function was implemented till late 1990 using an enhanced version of the program described in [4], incorporating all three abovementioned methods. This version, known as TC version 2.3, was developed by Alford and Havas in 1980. Features beyond those in [4] include a more space efficient coincidence procedure (along the lines described here) and somewhat improved coset definition strategies for the Felsch method (precursors to methods detailed here).

Since November 1990 a new version of coset enumeration, based on the principles described in this paper, has been incorporated in Cayley, and is now also available in a standalone form. HLT- and lookahead-type programs are considered together in this paper (generally under the name HLT-type), but special attention is given to Felsch-type methods. In broad detail the new version runs significantly faster than previous versions and, in some cases for Felsch-type enumerations, makes substantially fewer redundant coset definitions. Typically HLT-type enumerations run in about 60% of the time taken by the previous version while Felsch-type enumerations run in about 75%, except where there are significantly reduced numbers of redundant cosets in which case there are proportionately greater time savings.

---

[*]Supported in part by a Special Project Grant from the University of Queensland. This paper has benefited from detailed commentary on an earlier version by Charles Sims. Much of the programming and experimental work was done by Jin Xian Lian. This report appeared in *ISSAC'91* (Proceedings of the 1991 International Symposium on Symbolic and Algebraic Computation, Stephen M. Watt, editor), ACM Press, New York, 1991, pp. 191–199. Permission to copy without fee is granted by ACM provided that the copies are not made or distributed for direct commercial advantage.

The computational complexity of coset enumeration is not well understood. Even for a given strategy there is no computable bound, in terms of length of input and a hypothetical index, to the number of cosets which need to be defined in the coset enumeration process to obtain a complete table. (The existence of such a bound would violate unsolvability results for problems involving finitely presented groups.) Further, Sims [15] has proved that there is no polynomial bound in terms of maximum number of cosets for the number of coset tables which can derived by simple coset table operations like those used in coset enumeration programs. This more practical result indicates that the running time of coset enumeration procedures in terms of space available may be unpleasant. In this paper analysis of performance is given in terms of space and time requirements for selected examples chosen for illustrative purposes, including examples tabulated in [4]. Presentations for groups identified only by name in this paper are those given in [4].

The performance of the new program is demonstrated in enumerations of the 1140000 cosets of the subgroup generated by $\langle x, b, c, d, e, f, g \rangle$ in the Harada-Norton simple group, HN, presented by
$\langle x, a, b, c, d, e, f, g \mid x^2, a^2, b^2, c^2, d^2, e^2, f^2, g^2, (xa)^2,$
$(xg)^2, (bd)^2, (be)^2, (bf)^2, (bg)^2, (ce)^2, (cf)^2, (cg)^2, (df)^2,$
$(dg)^2, (eg)^2, (bc)^3, (cd)^3, (de)^3, (ef)^3, (fg)^3, (xb)^4, (ag)^5,$
$[a, edcb], [a, f]dcbdcd, [b, xcdefx], [cdef, xcdefx]\rangle$. These enumerations can be done readily in less than 40 cpu (and elapsed) minutes on a Solbourne 5/602 (about 22 mips) with 64 megabytes of real memory. Some details on such enumerations are given in §5 where they are used as test examples. The enumeration of the 55104 cosets of the subgroup $\langle b, c \rangle$ in Men(5), reported to have taken 79.5 hours on an IBM 360/50 in [4], is now readily done on a Sun 3/80 in about 3 minutes by both HLT- and Felsch-type methods. Of course the IBM 360/50 enumeration involved enormous amounts of paging, but it is now done in real memory.

## 2  Space considerations

In spite of the availability of much larger memories, space is still an important factor in coset enumeration programs. There are requests to compute indexes well into the millions, so efficient storing of coset tables is mandatory.

A significant advance in this direction is due to Beetham [2]. Coincidence processing techniques based on his methods are incorporated into modern coset enumeration programs. This means that the extra two columns used in [4] are no longer required.

The end result of these space saving techniques is that the coset table space required to enumerate and store $n$ cosets for a $g$ generator group which includes $i$ involutions is $(2g - i) \times n$ words rather than the previously required $(2g - i + 2) \times n$. Since $n$ can exceed a million in practical examples this saving can be very worthwhile. Observe that for a two generator group (without involutory generators) an extra 50% more cosets can be defined in a given amount of coset table space.

There are some costs associated with use of this space saving technique, costs worth paying. First, the coincidence routine is somewhat more complicated to code. Second, there is no longer easy availability of an active row list, which leads to a preference for compaction based storage reuse methods rather than free list based methods. This is exemplified in §5 in some large enumerations.

Atkinson [1] has pointed out that inverse columns may be dispensed with for any generators whose orders are known, not just for involutory generators. Thus, if the orders of $k$ out of $g$ generators are known, then enumerations can be carried out with $(2g - k)$ rather than $(2g - i)$ columns. However, for non-involutory generators, this is at the expense of considerable execution time and has not been utilized in the programs considered here.

## 3  Coset definition strategies

The order in which cosets are defined is explicitly prescribed by the order in which (the subgroup generators and) the group relators are processed in HLT-type procedures. In standard Felsch-type procedures coset definition is normally from the left/top of the coset table towards the right/bottom — that is, in order row by row. In fact a procedure needs to follow a method like this to some extent for the proofs that coset enumeration eventually terminates in the case of finite index (see [13]).

However this still leaves us with plenty of freedom in definition strategies, freedom which can be used to great advantage in Felsch-type methods. The work described here was embarked upon at the instigation of John Leech, who suggested that giving computer implementations more attributes of human enumerators (who define cosets intelligently rather than by rote) might be beneficial. Initial investigations began in 1979 leading to TC version 2.3 in 1980. More recent work leading to current strategies was undertaken in 1990.

Though it is not strictly necessary, Felsch-type programs generally start off by ensuring that each of the given subgroup generators forms a cycle at coset 1 before embarking on further enumeration. This is done by applying coset 1 to each subgroup generator, HLT-style, and it simplifies the ensuing logic since subgroup generators can be ignored thereafter.

In a complete coset table each group relator will eventually form a cycle at each coset, in particular at coset 1. It has been found to be a good general rule to use the given group relators as subgroup generators, that is, also form cycles at coset 1 for each group relator at the start of an enumeration. This ensures the early definition of some useful cosets.

Thus, as a default strategy (but under user control), the new version automatically treats all of the group relators as subgroup generators. Over a wide range of examples this has produced generally beneficial results in terms of maximum and total cosets. Thus, where there were few redundant cosets, the number of redundant cosets remained small. In pathological cases the tendency is for this strategy either to reduce the number of redundant cosets or to leave the number little changed.

A spectacularly good example of this behaviour is as follows. In [4] it was pointed out that some enumerations in the Macdonald groups $G(\alpha, \beta)$ provide examples where Felsch-type methods define significantly more cosets than HLT-type methods. Thus the maximum/total numbers of cosets defined for the index 40 subgroup $\langle [a,b], [b,a^{-1}], [a^{-1},b^{-1}], [b^{-1},a] \rangle$ of $G(3, 21)$ are 84/91 for HLT and 16063/16067 for Felsch when the relators are not used as subgroup generators. Use of the relators as additional subgroup generators reduces the Felsch figures to 56/59. (These figures do not incorporate the preferred definition strategies described later. Using the default preferred definition strategies the Felsch figures are 10586/10591 without relators as subgroup generators, 40/43 with.)

Over a substantial range of test examples, including all of the enumerations in [4], only one has been found in which incorporating relators as subgroup generators has a substantially bad effect. Thus the enumeration of the 100 cosets of $J_2$ over $\langle a, b, b^{ca^{-1}c} \rangle$ has maximum and total coset statistics go from 1305/1315 (without relators as subgroup generators) to 2617/2644 (with). Again these figures are without use of preferred definition strategies. The preferred definition strategies save the day. Statistics for this enumeration, with relators as subgroup generators and using the default preferred definition strategy, are 563/581, only slightly worse than the 511/528 achieved with default preferred definition strategy and no relators in the subgroup. Best user selectable strategy with the current program gives 366/381. In this case the seven shortest relators are included in the subgroup, the preferred definition list is a stack dropping early elements, and the fill fraction (defined later) is 1/4 to 1/6.

A side effect of the incorporation of the group relators as subgroup generators in Felsch-type methods is that Felsch-type methods become susceptible to the ordering of group relators and to cyclic permutation of group relators (in the same way as HLT-type was shown to be in [4].) Thus enumerations of the index 625 subgroup $\langle b, b^a \rangle$ in optimal presentations for the class 4 quotient of B(2,5) (defined in [8]) show statistics variations from 2117/2218 to 7475/7483, depending on which permutations of the two long commutator relations are used. Simply interchanging the order of these two relations shows variations of as much as from 4567/4569 to 5663/5666. The general guidelines given in [4] for HLT (put relators in increasing order of length and permute them so that they include long common subsequences) seem appropriate for Felsch-type methods now.

The key to performance of coset enumeration procedures is good selection of the next coset to be defined. Thus Leech has simplified a number of coset enumerations by removing useless cosets needlessly defined by computer implementations. In [10] he showed how readily 401 definitions can be removed from an enumeration in $G^{3,7,17}$, reducing total cosets from 1544 to 1143. Then, in [11], he considered enumerations in the Fibonacci group $F(2, 7)$. This group is studied in detail in [6], where a "best" enumeration with a total of 327 cosets defined (using programs of the day) was used to derive formally a proof of the structure of the group. Leech, by hand, reduced the total first to 129, then by pruning to 55. (This has subsequently been reduced to 53 by Edeson [5].)

Hand enumerators intelligently choose which coset should be defined next, based on the value of each potential definition. In particular, definitions which close relator cycles (or at least shorten gaps in cycles) are favored. A definition which actually closes a relator cycle immediately yields twice as many coset table entries as other definitions. Human enumerators can observe this at no extra cost to the enumeration process. It is shown here how this can also be done at little cost in computer implementations.

In early experiments in 1979 a minimum length gap or the first gap of length one which was detected in a relator cycle during deduction processing was noted. (This takes negligible extra time or space.) Then, when it came time to define a new coset, it was defined to fill this gap, if the gap still existed. (Subsequent deduction processing quite often filled the gap.) This led to some improvements in performance on pathological enumerations, improvements reported by the author in an unpublished lecture at the Symposium on Computational Group Theory in Durham in 1982. In general terms those improvements have now been overtaken.

In 1990 it became clear that all too often this one preferred definition, stored for use, was actually filled before definition time. It was therefore decided to store a "list" of preferred definitions. This still had to be done quickly and without consuming too much space, which led to gaps only of length one being retained. It immediately raised a question. How long might such a list of gaps of length one grow during a deduction processing phase in an enumeration?

In order to construct the list quickly it was decided not to check for repetitions (too time consuming). It was discovered that the list could become very long. Thus, during an enumeration of the 18468 cosets of $\langle a, b, c, s \rangle$ in $J_3^*$ (maximum 18468, total 18557) a list of preferred definitions of length 200614 (including repetitions) occurred. During an enumeration of the 29 cosets of the trivial subgroup in $F(2,7)$ (maximum 30964, total 30966) a list of preferred definitions of length 120745 (including repetitions) occurred.

Since the list of preferred definitions could become so long it was decided that it should be truncated in a sensible way. The idea is to keep enough to ensure that a preferred definition, assuming that there is one, is most likely to be made.

There are more questions about this list. What data structure should be used? Should it be a stack? Or a queue? Since it is being truncated, should items be discarded from the front? Or from the back? (Whether more complicated stuctures, which may include additional information, would provide better performance has not yet been yet been investigated.)

A naive expectation might suggest that "stack dropping earliest elements" is a good structure. In this case candidates for definition are tried using the latest found preferred position. However the details of deduction processing suggest that this might not be such a good idea. This is because relators are always processed in the same order, which could lead to preferred definitions from mainly later relators being at the top of the stack and earlier relators being relatively ignored.

Experiments were done with varying length lists and each of four data structures for the list: stack with earliest elements dropped; stack with latest elements dropped; queue with earliest elements dropped; queue with latest elements dropped. It was found that a list length of 200 was adequate in practice for a preferred definition to be usually available if any at all existed (though list length dependent on presentation length has intuitive appeal). With length 200 a comprehensive study was made of enumerations using each of the four mentioned data stuctures for the preferred definition list.

It is appropriate at this stage to mention a consequence of using such a preferred definition list. This kind of list used carelessly could lead to violation of the conditions required for guaranteed termination of a coset enumeration procedure in the case of finite index (mentioned in the first paragraph of this section.) To avoid such a violation all implementations described here include a mechanism to ensure that at least a certain initial fraction of the coset table (in order row by row) is filled before a preferred definition is made. Otherwise a standard definition is made. It is easy to construct examples where omission of a control of this type does indeed lead to infinite looping. This "fill fraction" is a user selectable parameter in these implementations.

The practical study of enumerations included a comprehensive variation of number of relators in the subgroup, data structure used for preferred definition list, and fill fraction. As a result the following parameters were chosen as default: all relators are used as subgroup generators; the preferred definition list is a queue dropping earliest elements; the fill fraction is chosen to be $1/\lfloor 5(c+2)/4 \rfloor$ (where there are $c$ columns in the coset table). Analyses giving practical grounds for the selection of these parameters are presented in [7]. (The fill fraction is chosen to be the reciprocal of an integer to enable implementation via fast code, without unduly affecting flexibility of choice. In practice the inverse of the fill fraction, the fill factor which is an integer, is used in the code.)

With the default strategy, over a substantial set of examples, it was found that the percentage of cosets filled from the preferred definition list (as against the standard order, row by row) varied greatly from enumeration to enumeration. Almost the whole range was represented, from 100% preferred definitions down to 10%. At most a small percentage (less than 10%) of definitions were forced away from preferred ones by the fill fraction.

As a challenge to computer enumerators Leech [11] set the target of 129 cosets for an enumeration in $F(2,7)$. The best achieved by the new Felsch-type program, with strategies selected by hand, not default, is a total of 169 cosets. We are getting there. (Strategy details are: no relators in subgroup; preferred definition list a length 200 queue dropping earliest elements; fill fraction 1/8.) This is the one example considered in this paper where (implicit) ordering of columns is taken into account (see [11]). The default strategy does not star: its best (taking simple cycling of columns by change of subgroup generator into account) is a maximum and total of 362 cosets.

# 4 Improving performance

Since 1973 technological advances in computers have dramatically changed the working environment for typical coset enumeration programs. Major changes apart from speed are the availability of much larger physical memories and the availability of even more substantial virtual memories. Care needs to be taken to gain proper advantage from these improvements. These technological changes are of course universal, not simply related to coset enumeration, and techniques used elsewhere apply.

It is appropriate to note at this stage that the 1973 implementation and the 1980 implementation were both written in FORTRAN (as was Cayley in 1980). The current implementation is written in C (as Cayley now is). This does not imply fundamental differences, rather differences in detail, which are mainly ignored here.

In 1973 space was the most precious resource for coset enumeration programs. This led to particularly space economical representations for all components of the data structures, including the presentation itself. Space savings which involved keeping relators in a base relator plus exponent form cost some extra execution time in terms of additional loops. With the multi-megabyte memories now available it is more advantageous to expand relators fully to speed processing. This, in effect, is a kind of loop-unrolling, common in execution speed-up techniques (see, for example, [9]). The approximate time saving for HLT-type implementations provided by this is about 10% overall.

The above applies to both HLT- and Felsch-type methods. An additional loop enhancement is applicable in Felsch, which involves expanding each relator out three times to avoid having to check for end of relator during deduction processing. This reduces the amount of code executed for each loop over each relator. For Felsch-type enumerations the combination of these loop enhancements leads to about 5% time savings overall.

A more significant saving comes from an analysis of the logical versus physical structure of the coset table. The coset table is logically a 2-dimensional array, normally stored with rows representing cosets and columns representing group generators or inverses. Entries in the coset table are coset numbers, frequently used for subsequent access to the coset table. The original 1973 and 1980 implementations used 2-dimensional arrays with FORTRAN facilities for array access. (Already in the FORTRAN versions the array was physically stored in logical row major order rather than column major order to reduce paging by enhancing locality of reference.) When ported to C for Cayley, analogous access was made via macro definitions for the C version of the 1980 implementation.

Let the coset table array be called **CT** and assume that there are $c$ columns. Thus access to **CT**$[i,j]$ (the action of the generator or inverse corresponding to column $j$ on coset number $i$) required the computation of $(i-1) \times c + j$ (assuming 1-base array subscripting). The replacement of the coset number $i$ in the coset table by $(i-1) \times c + 1$ (the $+1$ in this expression leaves coset 1 still represented by 1 and avoids zero, which means undefined, so that coding changes are relatively easy) leads to big savings in avoided multiplications by $c$ for coset table addressing. For HLT-type enumerations this device leads to execution time savings of about 25% while for Felsch-type the savings are about 15%. (The improvement for HLT is greater than that for Felsch because Felsch makes relatively fewer coset table accesses and does more work per access.)

There is a time saving observation which applies to coincidence processing. When total collapse occurs, to just one coset, it can happen that the result is known well before the full coincidence processing is complete. That is, only coset 1 is active (with all entries in row 1 defined) but there is a significant queue of deduced coincidences still waiting to be processed. In this case these coincidences can be discarded. A quick check for this situation is included in the implementations discussed here. Use of this check reduces the (Sun 3/80) execution time for an index 1 HLT-type enumeration of $G^{3,7,17}|E$ (maximum 122661, total 196972) by about 10% and the corresponding Felsch-type enumeration (maximum 57284, total 57592) by about 5% Even though, in the case of Felsch-type, the standard deduction queue is also discarded rather than processed, the time saving depends mainly on the size of collapse, not the enumeration type. (It is also possible to consider modifying the deduction queue in Felsch-type methods after any substantial collapse, but this is not investigated here.)

Further details about these programming modifications and some others are given in [7]. The end result is that the 1990 HLT-type implementations run about 40% faster than the 1980 version, while for Felsch the saving is about 25%, taking all the technical changes into account.

# 5 Virtual Memory

Paging may influence the choice of method. HLT-type procedures certainly define more cosets in a given amount of cpu time than Felsch-type procedures. However they generally lead to more cosets in total and more memory accesses. When large enumerations are anticipated, especially if the total number of cosets to be defined could be significantly larger than the number that will be resident in real memory, Felsch-type methods may well minimize paging. Generally speaking, HLT-type methods perform most quickly for "small" enumerations, but Felsch-type methods are often better for "large" or pathological enumerations. Felsch-type methods keep accesses more local than HLT-type. Perhaps the best solution is to avoid paging if possible by restraining table sizes to resident memory or thereabouts. A study of index 1140000 enumerations in HN (defined in the introduction) is instructive.

The performance figures should be read with the usual caveats about variability of such measurements, especially their dependence on system load. Most enumerations were run on a Solbourne 5/602 with 64 megabytes of real memory under similar load conditions (very little other activity). This meant that all of the program and data (about 40 megabytes) could usually fit into real memory with a 10 million word table size, but there was a deficit of 15 to 20 megabytes when a 15 million word table size was used (program plus data about 60 megabytes.)

On the Solbourne, Felsch-type enumerations behaved pretty well and consistently, regardless of virtual memory allocated and storage reuse method. Total cosets defined was about 1.5 million, which could fit into about 12 million words without coset reuse. Thus, with table sizes of both 10 and 15 million words and with either compaction or free list storage reuse, cpu time and elapsed time ranged from about 2300 to 2600 seconds. There were not many page faults (at most a few thousand), even with the 15 million word table size enumerations.

For these tests a complete (see [4]) version of lookahead was used. This is more time-consuming than methods proposed by Leon in [12] and emphasises the impact of lookahead on enumeration performance. With complete lookahead, a 10 million word table size and compaction, similar timings to Felsch were achieved. In total about 3.7 million cosets were defined. About 15% of the time was spent in the definition phase, 10% in compaction, and 75% in lookahead. However, when a free list was used it was quite different. Because of the new coincidence routine there was no active row list to enable coset application to relators in order of definition. Instead, cosets were applied to relators in their order in virtual memory. This led to a less efficient sequence of definitions increasing total cosets somewhat, to about 3.9 million. But it had a much worse impact on each lookahead pass. Since the cosets were processed in a less effective order, fewer coincidences were found per pass. The upshot of this was a substantial increase in the number of time-consuming lookaheads, from eight with compaction to fifteen with lookahead. There were still few page faults, but cpu and elapsed time nearly doubled, to about 3800 seconds. About 5% of the time was spent in the definition phase, and about 95% in lookahead.

Given a 15 million word table size, lookahead behaved poorly. With compaction, a total of about 4.8 million cosets were defined. Cpu times of around 4000 seconds were achieved, but elapsed times were around 70000 seconds. There were about 2 million page faults and more than half of the cpu time was system time. Thus about 1700 cpu seconds were user time and 2300 system time, essentially all paging. There were only four lookaheads, but they consumed most of the time. The definition phase took about 20% of the user time, compaction 10% and lookahead 70%. The definition and compaction phases took about 10% each of both system and elapsed times, with the remaining 80% taken by compaction. With a free list, this enumeration did not complete in some 160000 elapsed seconds which were available, at which stage over 1000 seconds of user time, over 4000 seconds of system time, and over 4 million page faults had been logged. In this case the definition phase had taken about 10% of the user time and about 5% of the system and elapsed times, with lookahead taking the remainder.

Disastrous behaviour was observed on a Sun SparcStation 2 with 40 megabytes of real memory. In spite of the availability of plenty of real memory the so-called 'PMEG thrashing problem' meant that a lookahead run with a 10 million word table size was restricted to about 1 megabyte of physical memory. It took 110 cpu **hours**, 35 elapsed **days**, and logged up some 60 million page faults. More than 107 hours of the cpu time was system time. In the same context Felsch was much better than lookahead, but still much worse than when given ample memory. Felsch took about eight cpu hours (seven of them system time), 25 elapsed hours, and half a million page faults.

## 6  Performance statistics

In this section an indication is given of the performance of the new program in a tabulation of maximum and total cosets defined and times taken for the pathological enumerations in Table 4 of [4], in the same order as in that paper. More detailed statistics are given in [7]. Here a comparison is made of the performance of the following: (A) a strategy equivalent to the 1973 Felsch program; (B) a strategy which simply adds all group relators as subgroup generators; (C) the default strategy, which does this and also uses a length 200 queue dropping earliest elements as a preferred definition list, with fill fraction of $1/\lfloor 5(c+2)/4 \rfloor$; (D) the best achieved with the same data structure as C, by manually setting the number of relators to be included in the subgroup and the fill fraction (note that a fill fraction of 1 disables use of the preferred definition list); (E) the best achieved across the four tested preferred definition list data structures by manually setting parameters as in D; and (F) the worst achieved that way. (This is purely a comparison of coset definition strategies since all of the programming suggestions of the previous section are incorporated into the coset enumeration implementation used.) Times taken are measured on a Sun 3/80 in cpu seconds. Note that as far as timing is concerned, for a given enumeration the time taken is basically proportional to the total cosets defined. Any one of these strategies can be selected by choosing parameters in the current version of the program. To find the best and worst requires running through all possible strategies. In practice relators were selected for incorporation into the subgroup in their given order simply by count and the lowest fill fraction used was 1/20.

These results are given in a double table on the next page, with the 13 enumerations identified by a number and the index. The top half shows how much improvement has been achieved by the new default strategy over previous methods. Then the bottom half shows how much more can be achieved reasonably readily, but also shows how badly a poorly chosen strategy can perform. As yet no generally good method of choosing the best strategy *a priori* rather than *a postiori* is at hand. The default strategy is the current *a priori* attempt at picking a good strategy.

It must be emphasised that this is a very restricted set of examples, chosen to illustrate behaviour because of its role for that purpose in [4]. The first thing to note is how large the differences can be between the best and the worst strategies. Since coset enumeration time depends mainly on total cosets defined, the ratio between totals of cosets for a given enumeration can be used as a performance measure for different methods of enumeration. Observe that the while in some cases (2, 4, 5, 6 and 7) the ratio between E (best) and F (worst) is less than two, in the other cases it is greater, exceeding twelve in case 3. In six of the eight cases where that factor exceeds two, method A is within a factor of two of the worst achieved. Only two times out of eight is Method A within a factor of two of the best achieved. Method C, the current default, is never within a factor of two of the worst for these eight enumerations, and is within a factor of two of the best seven times out of eight. Method C performs worse than A or B only once in this set of examples, by a factor not too far from one.

In [4] it is stated that "it appears that the lookahead algorithm can usually perform an enumeration in about the same space as required by the Felsch algorithm and occasionally in significantly less space." This is no longer true for the new strategy.

Significantly less space for lookahead arose for enumerations in Macdonald groups. This is remedied for Felsch-type methods by the use of relators as additional subgroup generators. Another example of the same kind of behaviour is given by Sims in [15]. He considers a presentation for the cyclic group of order $n^2 - 1$, namely $\langle x, y \mid x^n = y, y^n = x \rangle$. He points out that, for enumerations over the trivial subgroup, some versions of Felsch-type procedures lead to the definition of a total number of cosets which is exponential in $n$, while some versions of HLT define a total which is $O(n^3)$. Without relators included as subgroup generators, the new version of Felsch displays this exponential behaviour. However once the relators are included as subgroup generators (as indeed is done in the new default strategy), Felsch-type methods can complete the enumerations without defining any redundant cosets at all.

Compare best lookahead figures ($\hat{M}_L$ and $\hat{T}_L$ in terms of [4]) with the maximums and totals for the various Felsch strategies. There are four cases where $\hat{T}_L$ exceeds the Method C total by a factor of more than 2, six cases where $\hat{T}_L$ is somewhat greater (factor 1.1 to 2), two cases of approximate equality, and only one case where $\hat{T}_L$ is somewhat less than the Method C total. (The case where $\hat{T}_L$ is better, enumeration 12, is the Macdonald group enumeration $G(2,6)|E$, where $\hat{T}_L$ is 4582 and the Method C total is 7925.) When, perhaps more fairly, the best lookahead total is compared with the best Felsch-type method, E, six enumerations have $\hat{T}_L$ exceeding the Felsch total by a factor more than 2.

It should be emphasized that a number of relevant matters have not been considered in detail in this account. For the table, all enumerations used the subgroup generators and group relators explicitly as given, in the order given. Possible effects due to different subgroup generating sets, different relators, ordering of relators, inversion of relators, different presentations and order of columns have not been investigated. In general coset enumeration programs should massage both the subgroup generators and the relators before embarking upon serious enumerations. Indeed sorting both relators and subgroup generators by length is the (overridable) default strategy for the new implementation, but is not used in the tables presented here.

How well can a given enumeration be done? This is clearly a hard question. Witness the effort taken to reduce an $F(2,7)$ enumeration to a total of 53 cosets ([6], [11], [5]). And we are nowhere near a proof that 53 is best, should that indeed be the case. A combinatorial search for the best sequence of definitions, even knowing that there is bound of 53, has not been attempted. However it would probably be no mean task. (The first definition can go in any one of 12 places, the second in any one of about 20 places, ...)

How good are the preferred definition list methods described here? From the table we can see that for one example there is much improvement possible. For $G(2,6)|E$, $\hat{T}_L$ is 4582 while the best achieved by the variations tested here is 7893 (default, but fill fraction 1/3). Sims [14] has described situations where the Knuth-Bendix procedure outperforms coset enumeration. In [15] he says "where the final answer is small but long words are involved at intermediate stages of the verification, the use of the Knuth-Bendix procedure should be considered seriously". $G(2,6)$ is a case in point.

Table 1: Comparative Behaviour of Coset Enumeration Strategies

| Enum. | Index | $\hat{M}_L$ | A (1973) | | | B (relators) | | | C (default) | | |
|---|---|---|---|---|---|---|---|---|---|---|---|
| | | | Max. | Total | Time | Max. | Total | Time | Max. | Total | Time |
| 1 | 1 | 695 | 588 | 588 | 0.45 | 420 | 424 | 0.35 | 166 | 179 | 0.18 |
| 2 | 1 | 224 | 229 | 229 | 0.20 | 232 | 232 | 0.20 | 227 | 227 | 0.18 |
| 3 | 1 | 1381 | 1471 | 1471 | 2.08 | 1341 | 1342 | 1.85 | 1810 | 1837 | 2.52 |
| 4 | 660 | 661 | 981 | 1017 | 1.47 | 760 | 827 | 1.27 | 660 | 743 | 1.17 |
| 5 | 1092 | 2286 | 1775 | 1812 | 1.88 | 1497 | 1582 | 1.62 | 1251 | 1324 | 1.45 |
| 6 | 720 | 721 | 980 | 1223 | 2.60 | 1081 | 1385 | 2.97 | 720 | 837 | 2.17 |
| 7 | 448 | 1241 | 1302 | 1306 | 1.05 | 1287 | 1290 | 1.10 | 861 | 877 | 0.75 |
| 8 | 240 | 4553 | 4439 | 4740 | 7.55 | 4901 | 5234 | 8.45 | 4263 | 4333 | 7.97 |
| 9 | 120 | 2189 | 1638 | 1660 | 1.53 | 1374 | 1379 | 1.27 | 653 | 660 | 0.68 |
| 10 | 21504 | 69990 | 76308 | 76308 | 115.42 | 72703 | 72711 | 110.38 | 38250 | 39580 | 64.77 |
| 11 | 3 | 2973 | 6812 | 6864 | 6.05 | 5535 | 5562 | 4.88 | 3267 | 3279 | 3.13 |
| 12 | 5 | 4194 | 19597 | 19627 | 17.78 | 19646 | 19817 | 17.90 | 7921 | 7925 | 7.30 |
| 13 | 16 | 29007 | 109538 | 110105 | 93.52 | 50999 | 51064 | 44.72 | 26031 | 26114 | 23.85 |
| Enum. | Index | $\hat{T}_L$ | D | | | E (best) | | | F (worst) | | |
| | | | Max. | Total | Time | Max. | Total | Time | Max. | Total | Time |
| 1 | 1 | 758 | 112 | 127 | 0.10 | 98 | 104 | 0.10 | 588 | 588 | 0.45 |
| 2 | 1 | 227 | 216 | 216 | 0.18 | 216 | 216 | 0.18 | 242 | 242 | 0.18 |
| 3 | 1 | 2315 | 986 | 1008 | 1.33 | 724 | 761 | 1.00 | 8836 | 9218 | 11.72 |
| 4 | 660 | 824 | 660 | 743 | 1.12 | 660 | 743 | 1.12 | 1009 | 1050 | 1.47 |
| 5 | 1092 | 2880 | 1251 | 1324 | 1.45 | 1221 | 1310 | 1.45 | 1836 | 1868 | 1.92 |
| 6 | 720 | 1349 | 720 | 830 | 2.12 | 720 | 724 | 1.92 | 1081 | 1385 | 3.13 |
| 7 | 448 | 1422 | 861 | 877 | 0.75 | 824 | 840 | 0.75 | 1302 | 1306 | 1.05 |
| 8 | 240 | 7158 | 2928 | 2996 | 5.98 | 2650 | 2750 | 5.52 | 9371 | 9465 | 16.72 |
| 9 | 120 | 2206 | 653 | 660 | 0.68 | 653 | 660 | 0.68 | 1819 | 1876 | 1.65 |
| 10 | 21504 | 161805 | 21504 | 24642 | 41.23 | 21504 | 23702 | 39.33 | 131776 | 138218 | 184.05 |
| 11 | 3 | 3255 | 3188 | 3193 | 2.88 | 3188 | 3193 | 2.88 | 6812 | 6864 | 6.05 |
| 12 | 5 | 4582 | 7889 | 7893 | 7.48 | 7889 | 7893 | 7.48 | 19646 | 19817 | 17.02 |
| 13 | 16 | 31993 | 25498 | 25516 | 23.42 | 25481 | 25496 | 22.35 | 109538 | 110105 | 93.52 |

Tests in 1979 investigated both more complicated strategies and simpler strategies. More complicated strategies aimed at ensuring that preferred definitions were sequentially selected from different relators, to ensure that consequences of all relators were incorporated into the table. This was more time consuming and seemed to yield little appreciable improvement. However one such attempt gave the best known figures for enumeration 7 (810/824), although not really much better than Method E. Prior to the introduction of the preferred definition list, experiments were done with just one saved preferred definition. In this case it was easy to save a definition which filled a minimum length gap rather than just a length 1 gap. This proved to be a very erratic strategy, sometimes producing bad enumerations from good with just minor changes to parameters like the fill fraction. However, this strategy sometimes performed exceptionally well. Best computer figures for enumerations 2 (157/160), 5 (1092/1158), and 9 (438/438) have been achieved with variants of this.

How much improvement can be achieved by setting preferred definition list parameters well? How badly can enumerations go? From the table we see a factor exceeding 12 for enumeration 3, between best (no relators in subgroup, preferred definition list a length 200 stack deleting latest elements, fill fraction 1/3) and worst (only first relator in subgroup, same preferred definition list, fill fraction 1/14). The largest factor so far observed using a reasonable variation of parameters for the methods described here is almost 58. Thus, the index 105 enumeration $\langle x, y, z \mid x^y x^{-3}, y^z y^{-2}, z^x z^{-4} \rangle$ over $\langle x \rangle$ had best statistics 1648/1652 and worst 95425/95786. In both cases a preferred definition list length 200 queue deleting early elements was used. Best figures were achieved with only the first relator in the subgroup and a wide range of fill fraction, 1/6 or less.

The fill fraction was the real key to the differences, with good figures being achieved regardless of other parameters once the fill fraction was small enough to stop interference with use of the preferred definition list. The worst figures came with all relators in the subgroup and fill fractions 1 to 1/3. In these cases the large fill fraction prevented any use of the preferred definition list. Note also that the default Method C performs quite reasonably, achieving statistics 2854/2859, using a fill fraction of 1/10.